\documentclass[12pt,a4paper,reqno]{amsart}
\usepackage{amsfonts,amsthm,amsmath,amssymb,amscd,amsbsy,euscript,textcomp,graphicx}
\usepackage{color}
\usepackage[utf8]{inputenc}
\usepackage{url}

\newtheorem{definition}{Definition}
\newtheorem{problem}{Problem}
\title[Maximizing color difference in metro maps]{Painting new lines: \\[5pt]
 Maximizing color difference in metro maps
}
\author[S.\,F.\,Griffioen and A.\,V.\,Kiselev]{Simone Griffioen and Arthemy Kiselev}
\email{s.f.griffioen@rug.nl}
\address{Johann Bernoulli institute for mathematics and computer science\\
PO Box 407\\
9700 AK Groningen  }
\date{\today}
\begin{document}

\begin{abstract}
Each metro line usually has its own color on the map. For obvious reasons, these colors should be maximally different. Suppose a new metro line is built. We now explain a strategy to choose the color or several colors for the new lines, making them as different as possible from both the old ones and each other. This is illustrated by using the Moscow metro map.
\end{abstract}
\maketitle
\section{Introduction}
\noindent%
\begin{itemize}
\item[{\color[rgb]{0,0.6,0.36} $\bullet$}] ``Have you seen the new terminal on the \emph{green}?''
\item[{\color[rgb]{0.95,0.07,0.17} $\bullet$}] ``The \emph{red} is overcrowded every morning at this time.''
\item[{\color[rgb]{1,0.66,0.3} $\bullet$}] ``I live at the very end of \emph{orange}, still OK with me: the escalator brings me right into the office.''
\item[{\color[rgb]{0.73,0.13,0.54} $\bullet$}] ``I have to change to \emph{purple} here. See you tomorrow!''
\end{itemize}

Sure, it is all about metro. The subway. The metropolitan underground railway system. Each line has its own color and each color is the name of a line. A problem arises when another line is built: the old lines must keep their colors whereas the new line should have a color which looks maximally different from the old ones. It gets even more complicated if there are two or more new lines in a subway. Not only the new colors must be unlike the old ones, but also they must differ from each other as much as possible.

\begin{figure}
  \centering
  \includegraphics[scale=0.2]{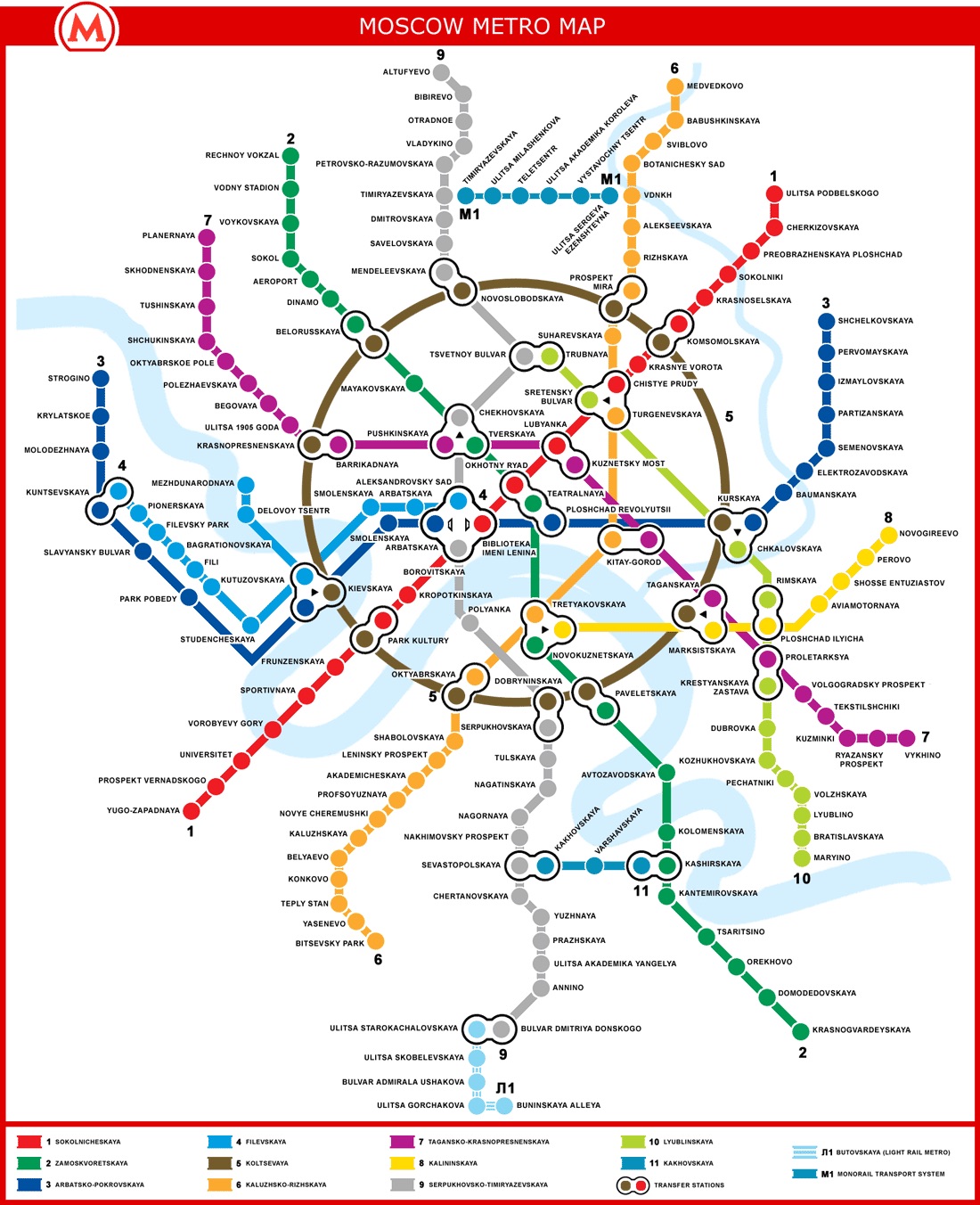}
  \caption{It is likely that you have already recognised the Moscow metro\protect\footnotemark\:here.} \label{moskoumetro}
\end{figure}

In this paper we discuss a mathematical strategy to choose these colors. By the way, what new colors would \emph{you} suggest for the real metro map in Figure \ref{moskoumetro}? For our purposes this map is ideal: not only there are 14 colors in it, but also there will be several new lines built in that subway in the near future (according to \cite{Reference8}). Let us take this map to illustrate our strategy --- yet the reader can apply the techniques to whathever other map he or she likes.

\section{Measuring the color difference}
\noindent%
To maximize color differences -- here, in metro maps -- we need a way to measure them and argue about colors objectively. For this we describe them as tuples of numbers by using some coordinate system (or \emph{color space}). For examble, the RGB color space is well-known: in it every color is encoded by its red, green, and blue components. Unfortunately, in the RGB space the Euclidean distance is not proportional to the visible distance that people perceive in between colors. That is, the RGB space is not \emph{perceptually uniform} --- which we would like.
\footnotetext{Opened in 1935 with 11 km of lines and 13 stations, the Moscow metro has 196 stations nowadays and the total route length of its lines is at least 327,5 km \cite{Reference7}. It is one of the world's most heavily used metro systems.}

\begin{figure}
  \centering
  \includegraphics[]{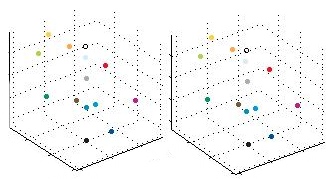}
  \caption{A stereographical image of the colors of the Moscow metro map at their location in the CIELAB space. The $L$-axis is the vertical axis, the $a$-axis runs from south-west to north-east and the $b$-axis from south-east to north-west. (\emph{To see a 3-dimensional picture, view the left picture with your left eye and the right picture with your right eye at the same time.})}\label{3dplaatje}. 
\end{figure}

In 1931, the Commission Internationale de l'Eclairage (CIE) did measurements on vision of hundreds of humans. From those data, several color spaces were created, to make communication about colors easy and unambiguous. In 1976 the CIE recommended the use of two approximately perceptually uniform color spaces and color-difference formulas. These spaces have become known by their officially recommended abbreviations CIELAB and CIELUV. In the former, a color is described by three coordinates: $L$ (lightness), $a$ (red-green scale) and $b$ (yellow-blue scale).\footnote{The L-scale runs from 0 to 100; the parameters $a$ and $b$ typically have values in the range -150 to 150. In this article we will round the coordinates of the colors to integers. The human eye can discriminate up to ten million colors, so it can distinguish the differences even finer than integer steps.} We will use this space to describe the colors of the Moscow metro map in terms of numbers: their CIELAB coordinates can be found in Table \ref{tabelmoscowpoints}; Figure \ref{3dplaatje} shows a stereographical plot of the colors of the Moscow metro in the CIELAB space.

\begin{figure}[h]
  \centering
  \includegraphics[scale=0.6]{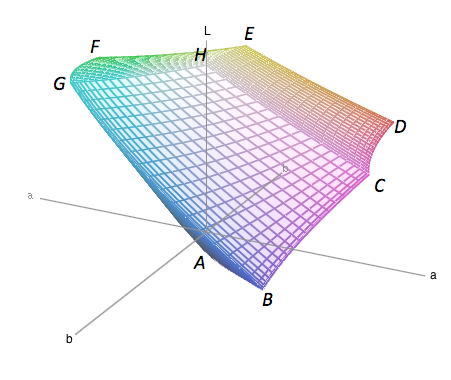}
  \caption{The gamut of Adobe 98 in the CIELAB space. We approximate this gamut by a polyhedron by taking the cornerpoints of this polyhedron close to the points $A,\ldots H$ in the figure.}\label{gamutAdobe}
\end{figure}

The CIELAB space itself is not bounded; however, based on what one wants to measure, only part of it matters in practice (for example, only the colors that a printer can print or only the colors the human eye can see). Such a subset of colors is called a color \emph{gamut}. The gamut of Adobe 98 in the CIELAB space can be seen in Figure \ref{gamutAdobe}. To keep the problem simple, we aproximate this gamut by a polyhedron. Because the CIELAB space was supposed to be almost perceptually uniform, the CIE 1976 color difference formula (CIE76) is given by the Euclidean distance
\begin{align*}
\Delta E_{ab}=\sqrt{(L_1-L_2)^2+( a_1-a_2)^2 + (b_1-b_2)^2}.
\end{align*}
In fact, the CIELAB space is not \emph{really} uniform. In particular, at high values of $a$ and $b$, the far too simple CIE76 formula values color differences too stronly compared to the experimental results on color perception. This is why several new color-difference formulas were proposed, of which the newest is the CIEDE2000 formula\footnote{In the CIEDE2000 distance some weighting coefficients have to be chosen. We use the ratio $k_L^{-1}:k_c^{-1}:k_h^{-1}=\frac{1}{2}:1:1$. This means that we consider the difference in \emph{chroma} $c=\sqrt{a^2+b^2}$ and \emph{hue} $h=$arctan$(b/a)$  to be more important than the difference in lightness $L$. For more information about this color difference formula, see \cite{Reference12}.} (or $\Delta E_{00}$). 

\begin{table}[]
\centering
\begin{tabular}{ r | l | r r r r }
  Line & Color & $L$ & $a$ & $b$& $\min{\Delta E_{00}}$ \\
  \hline
  1 & {\color[rgb]{0.95,0.07,0.17} $\bullet$} Red & 52 & 74 & 53 & 29,1 \\
  2 & {\color[rgb]{0,0.6,0.36} $\bullet$} Ocean green &56 & -45 & 26 & 21,7 \\
  3 & {\color[rgb]{0.2,0.33,0.57} $\bullet$} Cobalt & 35 & 7 & -43 & 19,4\\
  4 & {\color[rgb]{0,0.62,0.84} $\bullet$} Sky blue & 61 & -16 & -42 & 6,3\\
  5 & {\color[rgb]{0.46,0.36,0.22} $\bullet$} Olive brown & 41 & 6 & 27 & 20,4\\
  6 & {\color[rgb]{1,0.66,0.3} $\bullet$} Peach & 76 & 24 & 67 & 14,5\\
  7 & {\color[rgb]{0.73,0.13,0.54} $\bullet$} Pinkish purple & 43 & 64 & -24 & 28,3\\
  8 & {\color[rgb]{1,0.82,0.28} $\bullet$} Light mustard & 86 & 4 & 85 & 14,5\\
  9 & {\color[rgb]{0.67,0.68,0.69} $\bullet$} Light grey & 71 & 0 & -2 & 9,4\\
  10 & {\color[rgb]{0.69,0.82,0.33} $\bullet$} Greenish yellow & 80 & -26 & 68 & 17,8\\
  11 & {\color[rgb]{0,0.58,0.71} $\bullet$} Turquoise blue& 56 & -21 & -29 & 6,3\\
  Water & {\color[rgb]{0.82,0.94,0.98} $\bullet$} Pale blue& 93 & -8 & -9 & 11,3\\
  Background & $\circ$ White & 100 & 0 & 0 & 9,4\\
  Text & {\color[rgb]{0.14,0.12,0.13} $\bullet$} Black& 13 & 2 & 0 & 20,4\\
\end{tabular}
\caption{Coordinates of colors on the Moscow metro map. The names of the colors come from \cite{Reference10}. The average minimal distance is 16,3. Note that some of the minimal distnces are the same; this happens when the two points are the closest to each other.}
\label{tabelmoscowpoints}
\end{table}

\section{Mathematical formulation of the problem}
\noindent%
Given the CIELAB space equipped with the Euclidean metric -- or even better, the CIEDE2000 color difference formula -- and the color gamut $\Gamma$ (which we chose to be a polyhedron), we consider a set of existing colors $P=\{p_1,\ldots, p_k\}$ within this gamut, represented by their CIELAB coordinates (e.g., see Table \ref{tabelmoscowpoints}). If we want to add only one new line, this leads to the following maximin problem:
\begin{problem}[one new color]\label{onecolor}
How can we find an $x\in\Gamma$ such that
\begin{align*}
\min_{p_i\in P}{\|x-p_i\|}\to\max_{x\in \Gamma}\:?
\end{align*}
\end{problem}
\noindent%
That is, the minimal distance from the new color $x$  to the existing ones must be as large as possible, provided that the color $x$ stays inside or on the boundary of the gamut $\Gamma$. 

If we want to add several new lines, each having its own new color, the maximin problem becomes this.
\begin{problem}[$m$ new colors]\label{morecolors}
How can we find a collection of $m$ points $X=\{x_1,\ldots,x_m\}$ such that
\begin{align*}
\min_{\substack{x_i,x_j\in X\\i\neq j\\p_k\in P}}\{\|x_i-x_j\|,\|x_i-p_k\|\}\to\max_{X\subset \Gamma}\:?
\end{align*}
\end{problem}
\noindent%
That is, we maximize both the minimal distance in between  the new colors $x_1,\ldots,x_m$ and their minimal distance to the old colors. Alternatively, we could add the new colors $x_1,\ldots,x_m$ one by one, using a method for solving Problem \ref{onecolor}. However this does not always yield the optimal result.\footnote{Metro lines are built one after another, so choosing new colors one by one is quite logical. However, sometimes we know in advance that several new lines are being built.}

\section{The Voronoi diagram method}
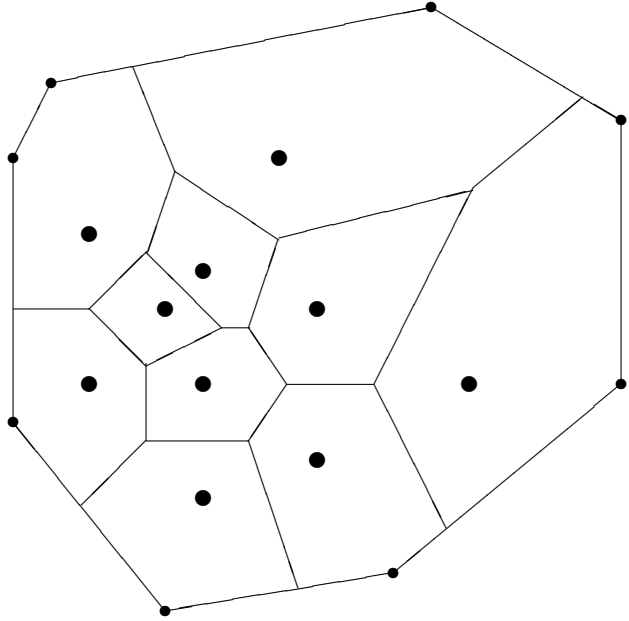
\begin{figure}
\setlength{\unitlength}{2mm} 
\begin{center}
\unitlength=1mm
\linethickness{0.4pt}
\begin{picture}(100,80)
\put(30,35){\circle*{2}}
\put(30,55){\circle*{2}}
\put(40,45){\circle*{2}}
\put(45,20){\circle*{2}}
\put(45,35){\circle*{2}}
\put(45,50){\circle*{2}}
\put(55,65){\circle*{2}}
\put(60,25){\circle*{2}}
\put(60,45){\circle*{2}}
\put(80,35){\circle*{2}}
\put(20,45){\line(1,0){10}}
\put(30,45){\line(1,1){7.6}}
\put(30,45){\line(1,-1){7.5}}
\put(37.5,27.5){\line(0,1){10.33}}
\put(37.5,27.5){\line(1,0){13.5}}
\put(37.5,27.5){\line(-1,-1){8.5}}
\put(56,35){\line(-2,-3){5}}
\put(42.5,40){\line(-2,-1){5}}
\put(42.5,40){\line(2,1){5}}
\put(47.5,42.5){\line(1,0){3.7}}
\put(47.5,42.5){\line(-1,1){10}}
\put(56,35){\line(1,0){11.5}}
\put(56,35){\line(-2,3){5}}
\put(51,42.5){\line(1,3){3.85}}
\put(41.3,63.25){\line(-2,5){5.5}}
\put(41.3,63.25){\line(-1,-3){3.6}}
\put(41.3,63.25){\line(3,-2){13.7}}
\put(55,54.4){\line(4,1){25.5}}
\put(67.5,35){\line(1,2){13}}
\put(80.5,61){\line(6,5){14.5}}
\put(67.5,35){\line(1,-2){9.5}}
\put(51,27.5){\line(1,-3){6.5}}
\put(20,30){\circle*{1.5}}
\put(20,65){\circle*{1.5}}
\put(25,75){\circle*{1.5}}
\put(75,85){\circle*{1.5}}
\put(100,70){\circle*{1.5}}
\put(100,35){\circle*{1.5}}
\put(70,10){\circle*{1.5}}
\put(40,5){\circle*{1.5}}
\put(20,30){\line(0,1){35}}
\put(20,65){\line(1,2){5}}
\put(25,75){\line(5,1){50}}
\put(75,85){\line(5,-3){25}}
\put(100,70){\line(0,-1){35}}
\put(100,35){\line(-6,-5){30}}
\put(70,10){\line(-6,-1){30}}
\put(40,5){\line(-4,5){20}}
\end{picture}
\end{center}
\caption{A 2-dimensional Voronoi diagram. Note that all edges of the Voronoi polygons are parts of the bisectors of lines connecting two points in the diagram.}
\label{2dvoronoi}
\end{figure}

\noindent%
Let us consider Problem \ref{onecolor}. Because the CIELAB space is continuous (at least, in principle), it is impossible to ``try all points'' in the gamut and find the one with the maximal minimal distance. Therefore it would be much easier to have only a finite number of candidate points to check. This can be achieved by using Voronoi diagrams.
\begin{definition}
Let $P=\{p_1,\ldots, p_k\}\subset\mathbb{R}^3$ be a finite set \textup{(}e.g., of old colors\textup{)} such that $p_i\neq p_j$ for $i\neq j$. The \emph{Voronoi polyhedron} of a point $p_i$ is the region \begin{align*}
V(p_i):=\{x\in \mathbb{R}^3:\: \|x-p_i\|\leq \|x-p_j\| \mbox{ for all } j\neq i\}.
\end{align*}
The set $\mathcal{V}(P)=\{V(p_i),\ldots, V(p_k)\}$ is the \textup{(}3-dimensional\textup{)} \emph{Voronoi diagram} generated by the set of points $P$. 
\end{definition}\label{voronoidiagram}
It is easily shown that the sought-for point $x$ maximizing the minimal distance in Problem \ref{onecolor} belongs to the intersection of
\begin{description}
\item[\textmd{either}] three Voronoi planes, that is, three faces of Voronoi polyhedra,
\item[\textmd{or}] two Voronoi planes and a boundary plane of the gamut,
\item[\textmd{or}] a Voronoi plane and two boundary planes,
\item[\textmd{or}] three boundary planes.
\end{description}
This means that as soon as these candidate points are calculated,\footnote{This can be done by elementary linear algebra, that is, by calculating bisectors and intersections of planes. In Figure \ref{2dvoronoi}, a 2-dimensional example of a Voronoi diagram can be found.} it only remains to compare their minimal distances to the existing points $p_i\in P$ and then, take the one whose distance is maximal; in fact, the resulting set is finite and reasonably small. This is the \emph{CIE76 Voronoi method} of solving Problem \ref{onecolor} with respect to the Euclidean distance. In Table \ref{optimalcolorsonebyone} we list the five new colors that we obtained, one by one, for the Moscow metro by using this method.

It is readily seen that the Voronoi method would provide the global optimal solution to maximin Problem \ref{onecolor} --- if there were a perceptually uniform color space. This is not yet the case nowadays, therefore let us study how better extimates can be obtained.

\begin{table}[h]
\centering
\begin{tabular}{ l | l | l | r r r}
  Color & $\min{\Delta E_{ab}}$ & $\min{\Delta E_{00}}$ & $L$ & $a$ & $b$ \\
  \hline
  {\color[rgb]{0,1,0} $\bullet$} Bright green & 114,4 & 22,0 & 83 & -138 & 91\\
  {\color[rgb]{0,0,1} $\bullet$} Blue& 87,0 & 19,8& 33 & 80 & -109\\
  {\color[rgb]{0.27,1,0.60} $\bullet$} Sea green & 67,7 & 13,1& 85 & -106 & 31\\
  {\color[rgb]{0,1,1} $\bullet$} Cyan & 65,3 & 21,9 & 87 & -78 & -21\\
  {\color[rgb]{0.85,0,1} $\bullet$} Bright pink& 56,1 & 14,4 & 60 & 100 & -64\\
\end{tabular}
\caption{Colors we find one by one for the first five new lines in the Moscow metro by using the CIE76 Voronoi method, that is, by solving Problem \ref{onecolor} with respect to the Euclidean distance. The minimal distances which we indicate here take in account the points that came earlier into the table. Note that the minimal $\Delta E_{00}$ distance is then not necessarily decreasing. }
\label{optimalcolorsonebyone}
\end{table} 

Because the CIEDE2000 distance is much more accurate than the Euclidean distance, we can improve our method by selecting the optimal point from the candidate points by comparing the $\Delta E_{00}$ distances instead of the Euclidean distances $\Delta E_{ab}$. This is the combined \emph{CIE76-CIEDE2000 Voronoi method}. Of course, this method does not solve Problem \ref{onecolor} with respect to $\Delta E_{00}$ alone, because the candidate points are found by using $\Delta E_{ab}$. However it is a fast method to obtain colors by using the CIEDE2000 distance. The five new colors obtained, one by one, for the Moscow metro by this method are given in Table \ref{optimalcolorsonebyoneciede2000}. We see that this gives us much better solutions with respect to the CIEDE2000 distance than the Voronoi CIE76 method: indeed, the $\Delta E_{00}$ distance drops at a much slower rate.

\begin{table}[h]
\centering
\begin{tabular}{ l | l | r r r}
  Color & $\min{\Delta E_{00}}$ & $L$ & $a$ & $b$ \\
  \hline
  {\color[rgb]{0,1,1} $\bullet$} Cyan & 23,9 & 87 & -78 & -21\\
  {\color[rgb]{0, 1,0} $\bullet$} Bright green& 22,0& 83 & -138 & 91\\
  {\color[rgb]{0, 0, 1} $\bullet$} Blue& 19,8 & 33 & 80 & -109 \\
  {\color[rgb]{0.9,0.79,1} $\bullet$} Pale lavender & 19,7 & 86 & 23 & -22\\
  {\color[rgb]{1,0.56,0.66} $\bullet$} Rose pink & 19,8 & 78 & 59 & 12\\
\end{tabular}
\caption{Colors we find one by one for the first five new lines using the combined CIE76-CIEDE2000 Voronoi method. Note that the distances here are the CIEDE2000 distances, and they drop at a much slower rate than in Table \ref{optimalcolorsonebyone}.\protect\footnotemark }
\label{optimalcolorsonebyoneciede2000}
\end{table}

\section{The simplex method}
\noindent%
We can also solve Problem \ref{morecolors} with respect to the Euclidean distance by using the simplex method. Let us consider here the case when the number of colors $m$ is two, which gives us the following function to maximize,
\begin{align*}
M(x,y)=\min_{\substack{p_i\in P\\x,y\in\Gamma}}\{\|x-p_i\|,\|y-p_i\|, \|x-y\|\}.
\end{align*}
We have $(x,y)\in \Gamma\times \Gamma\subset \mathbb{R}^6$. Let us rewrite our problem: first, put one extra parameter $x_0$ and the two vectors $x,y$ together in $\tilde{x}=(x_0,x,y)\in \mathbb{R}^7$, and now maximize the function $M_0(\tilde{x})=x_0$ over all $\tilde{x}\in \mathbb{R}\times \Gamma\times \Gamma$ satisfying the constraints
\begin{subequations}\label{nonsmooth1}
\begin{align}
\|x-p_i\|-x_0 &\geq0 \quad\mbox{ for each $p_i\in P$},\\
\|y-p_j\|-x_0&\geq0 \quad\mbox{ for each $p_j\in P$},\\
 \|x-y\|-x_0&\geq0.
\end{align}
\end{subequations}
Of course, more constraints on $(x,y)$ are used to delimit the set $\Gamma\times \Gamma$. Now we use the simplex method to maximize our function.

\footnotetext{The minimal CIEDE2000 distances do not always decrease each step; this is because the candidate points are selected by using the Euclidean distance, which does not necessarily give an optimum with respect to the CIEDE2000 distance. Therefore it is possible that, after adding a new color, one of the new candidate points is even a better solution than the point added before.}

The gamut's boundary gives the linear constraints, and the nonlinear constraint functions are equations \eqref{nonsmooth1}. Because the simplex method gives us only a local optimum, we apply it for a million random initial points $(x_0,x,y)\in\mathbb{R}^7$ (with $x_0\in[0,20]$ and $(x,y)\in\Gamma\times\Gamma$). These are local optima with respect to the $\Delta E_{ab}$ distance; let us now choose the one with maximal $\Delta E_{00}$ distance. This way, we obtain the two colors \[\mbox{Bright aqua  } {\color[rgb]{0.49,0.94,0.90} \bullet} \mbox{  (84, -65, -12)} \quad \mbox{ and }\quad \mbox{Baby pink  }{\color[rgb]{1,0.74,0.8} \bullet} \mbox{  (86, 35, 7)}.\] Their combined minimal CIEDE2000 distance to the old colors and each other is 22,7. This is more than the minimal distance which was obtained after two steps of the one by one CIE76-CIEDE2000 Voronoi method. So this gives us a better result for solving Problem \ref{morecolors}. However, this method is quite time-consuming (as $m$ increases, the simplex method will need even more time), so the much faster Voronoi method is a good alternative. 

\section{Conclusion and discussion}
\noindent%
We conclude that the Voronoi diagram method (and its improvement by using both the CIE76 and CIEDE2000 distance functions) for adding new colors one by one gives solutions which are comparable with the ones produced by the nonsmooth optimization method for finding two or more new colors at once. This is true for the map of the Moscow metro under study. Other solution schemes are of course possible, e.g. the \emph{Monte-Carlo} method: test points $x_i\in\Gamma$ for the new color(s) are generated within the gamut at random, and the new color -- or a tuple of new colors -- is the maximin optimum over a sufficiently large number of such attempts. The \emph{ballistic method} is realized as follows: first, suppose that each old color $p_i\in\Gamma$ and every new point $x_j$ in the sought-for tuple carries a positive electric charge (so that they repulse by Coulomb's law). To prevent an escape of new points from the gamut $\Gamma$ across its boundary, put a negligibly small but positive charge at the boundary's closest point to $x_j$ for every $j$. By starting from random initial data $x_j(0)\in\Gamma$, the configurations $\{x_j(t)\}$ of charges evolve inside the polyhedron $\Gamma$ which, by assuption, is filled with a viscous medium so that the new points' motion slows down; sooner or later they are captured near the local minima of the potential (see \cite{Reference20} for details on both methods).

Finding several new colors is still not all that one needs in order to paint the new lines on a metro map. Clearly, it is the configuration of transfers from the old to new lines that must further be taken into account --- to distribute the new colors between the new lines in such a way that at every transfer station, the colors of intersecting lines are maximally different.

\subsection*{Breaking news: an on-line poll}
\noindent%
The 16.8 km long Kozhukhovskaya line with eight stations on it is
scheduled to appear in the Moscow underground system in 2015-16,
officials say. The new line will connect the (south-)east with the
city centre.

Two colors, to freely choose from, were offered to pre-registered
``active citizens'' during an interactive poll, held from 20 October
till 10 November 2014. Namely, \emph{black} and \emph{pink} were
claimed to remain the only colors from the gamut that were not yet
used to paint lines on the Moscow metro map (see
Figure \ref{moskoumetro}). As many as 307,350 active citizens
took part in the on-line poll, media report; 71\% voted pink and 18\%
preferred black, whereas 7\% trusted the choice to experts (such as
designers, ethnographers, or psychologists). Finally, 4\% of the
respondents either made their own suggestions of the new color:
turquoise, olive, coral, etc. -- or even proposed to make that line on
the map speckled (see \cite{Reference22,Reference21}).

Permitting a large group of people (here, 307,350 pre-registered
active citizens) to choose the new color(s) from a given list of
pre-selected variants (here, two: black vs pink) might not always be
the optimal strategy to resolve an important issue, painting new
line(s) on a metro map. Perhaps, the decision should be trusted to the
Protector of tunnel builders and Best Friend of all cartographers;
then, the final choice is sometimes delegated to a child.
Alternatively, one could just ask a mathematician.


\begin{thebibliography}{10}
\bibitem{Reference20}Griffioen, S. (2014) Maximizing color difference in metro maps. \url{http://scripties.fwn.eldoc.ub.rug.nl/}.

\bibitem{Reference8} Laconte, P. Moscow: looking to the future - Mobility. \url{http://archive.ffue.org/PDF/Moscow-PL-02-rev.pdf}.

\bibitem{Reference12} Luo, M.R., Cui, G. and Rigg B. (2001) The development of the CIE 2000 colour-difference formula CIEDE2000. \emph{COLOR research and application 26}, 340--350.

\bibitem{Reference7} Moscow metro - official site. \url{http://engl.mosmetro.ru/}.

\bibitem{Reference50} Moscow.info. \url{http://http://www.moscow.info/}.

\bibitem{Reference10}Munroe, Randall. xkcd color survey results. \url{http://blog.xkcd.com/2010/05/03/color-survey-results/}

\bibitem{Reference5}Okabe, A., Boots, B., Sugihara, K., and Chiu, S.N. (2000)
 \emph{Spatial tessellations: Concepts and Applications of Voronoi Diagrams}. John Wiley and Sons Ltd, Chichester etc. 
 
\bibitem{Reference6}Robertson, A.R. (1990) Historical Development of CIE Recommended Color Difference Equations. \emph{Color research and application 15}, 167--170.

\bibitem{Reference16}Rockafellar, T. (1994) Nonsmooth optimization. \emph{Mathematical programming: State of the art}, 248--258.

\bibitem{Reference22} Ros Business Consulting. \url{http://www.rbc.ru/rbcfreenews/5464e5a3cbb20f3d459fa4a9}.

\bibitem{Reference21}Russian Info Agency ``Novosti".
\url{http://ria.ru/moscow/20141113/1033155060.html}.

\end{thebibliography}
\end{document}